\documentclass[11pt, a4paper]{article}
\usepackage{amsfonts,enumerate,amsmath,amsfonts,xspace,array,delarray,theorem,amssymb,epsfig}

\usepackage{graphicx,color}
\usepackage[latin1]{inputenc}%

\newenvironment{Myitemize}{ 
\begin{itemize}}{\end{itemize}}

\theoremstyle{break}\theorembodyfont{\it}
\theoremheaderfont{\normalfont\bfseries}

\newtheorem{theo}{Theorem}

\newtheorem{defi}[theo]{Definition}

\newtheorem{lem}[theo]{Lemma}

\newtheorem{prop}[theo]{Proposition}

\newtheorem{coro}[theo]{Corollary}

\newenvironment{proof}{\noindent{\bf Proof: }}
                {\leavevmode\unskip\nobreak\hskip2em plus1fill
                $\scriptstyle\square$\vskip\theorempostskipamount\par}

\newenvironment{rem}{\vskip\theorempostskipamount\par\noindent{\bf Remark: }}
{\leavevmode\unskip\nobreak\hskip2em plus1fill
\vskip\theorempostskipamount\par}

\let\lt=<
\let\gt=>
\def\freccia{{\longrightarrow}}
\def\R{{\mathbb R}}

\def\N{{\mathbb N}}
\def\Z{{\mathbb Z}}

\let\phi=\varphi
\let\eps=\varepsilon
\def\S{{\mathcal S}}

\def\be{\begin{equation}}
\def\ee{\end{equation}}
\def\beq{\begin{eqnarray*}}
\def\eeq{\end{eqnarray*}}
\def\e{{\rm e}}
\def\d{{\rm d}}

\begin{document}
\title{Propagation of Gevrey regularity for solutions of the Boltzmann equation for
Maxwellian molecules}
\date{}

\author
{
Laurent DESVILLETTES\\
{\small CMLA, ENS Cachan,}\\
{\small 61, avenue du Président Wilson, 94235 Cachan Cedex, France} \\
{\small E-mail:} \texttt{desville@cmla.ens-cachan.fr}\\
\\
 Giulia FURIOLI \\
{\small Dipartimento di Ingegneria Gestionale e dell'Informazione, Universit\`a
di Bergamo,}\\
{\small Viale Marconi 5, I--24044 Dalmine (BG), Italy} \\
{\small E-mail:} \texttt{gfurioli@unibg.it}\\ \ \\
 Elide TERRANEO\\
{\small Dipartimento di Matematica F.~Enriques, Universit\`a
degli studi di Milano,}\\
{\small Via Saldini 50 , I--20133 Milano, Italy} \\
{\small E-mail:} \texttt{terraneo@mat.unimi.it}
}
\maketitle

\begin{abstract}
We prove that Gevrey regularity is propagated by the Boltzmann
equation with Max\-wellian molecules, with or without angular
cut-off. The proof relies on the Wild expansion of the solution to
the equation and on the characterization of Gevrey regularity by
the Fourier transform.
\end{abstract}

\section{Introduction}
This paper deals with a propagation property for the solution of the following
Cauchy problem for the spatially homogenous Boltzmann equation for Maxwellian
molecules
\[
\left\{
\begin{aligned}
&\partial_t f(v,t)= Q(f,f)(v,t)\\
& f(v,0)= f_0(v).
\end{aligned}
\right .
\]
Here, $f(v,t) : \R^3 \times \R^+ \freccia \R$ is the probability
density of a gas which depends only on the velocity $v\in\R^3$ at
the time $t\geq 0$ and $Q$ is the quadratic Boltzmann collision
operator in the case of Maxwellian molecules:
\begin{equation}\label{Q}
Q(f,f)(v,t)=\int_{w\in\R^3}\int_{n\in S^2}
\left (f(v_*, t) f(w_*, t) -f(v,t)f(w,t)\right ) b\left(\frac{v-w}{|v-w|}\cdot n
\right ) \, \d n \, \d w.
\end{equation}
Due to the physical assumptions that the gas evolves through binary, elastic
collisions which are localized both in space and time, the relations between
the velocities $(v_*,w_*)$ of two particles before the collision and $(v, w)$
after it are the following~:
\[
\left\{
\begin{aligned}
& v_*= \frac {v+w} 2 + \frac {|v-w|}2 \, n,\\
& w_* = \frac {v+w} 2 - \frac {|v-w|}2 \, n,
\end{aligned}
\right .
\]
where $n$ is a vector in $S^2$, the unit sphere in $\R^3$, and parametrizes
all the possible pre-collisional veloci\-ties.

The collision kernel $b$, which is supposed to be nonnegative, is the function
which selects in which way the pre-collisional velocities contribute to produce
particles with velocity $v$ after the collision and is supposed (this is
precisely the assumption of Maxwellian molecules) to depend only on the cosine
of the deviation angle $\theta$, namely
\[
\cos\theta= \frac{v-w}{|v-w|}\cdot n.
\]
Finally, we will make the so-called non cut-off assumption, which means that
$b\notin L^1_{\rm loc} (]-1,1[)$ and, more precisely, we shall consider
\[
b(\cos \theta ) \sim \frac 1{ (1-\cos \theta)^{\frac 54}},\quad \theta \to 0.
\]

From a physical point of view, that means that the gas molecules
repel each other with a force proportional to the fifth power of
their distance and a great contribution to the integral collision
term is given by the grazing collisions ($\theta \sim 0$). The
assumption that the collision kernel $b$ is instead integrable on
$]-1,1[$ is called a cut-off assumption. For more information
about Boltzmann equation and its physical meaning, the reader can
consult for instance the review article by Villani \cite{vil}.

Due to the singularity of the collision kernel at the origin, the integral term
\eqref{Q} is not meaning\-ful if $f$ is not smooth and so it is convenient to
consider the weak form of the Boltzmann equation: for
$\phi \in C^\infty_c (\R^3)$,
\[
\begin{aligned}
& \int_{v\in \R^3} \partial_t f(v,t) \phi(v) \, \d v \\
&= \int_{v\in\R^3}
 Q(f,f)(v,t) \phi(v) \, \d v\\
& = \int_{v\in\R^3}\int_{w\in\R^3}
\int_{n\in S^2}
f(v,t) f(w,t)\left (\phi (v_*) - \phi (v)\right ) b\left(\frac{v-w}{|v-w|}\cdot n
\right) \, \d n \, \d w\, \d v,
\end{aligned}
\]
or even, with another point of simplification, in the Fourier variable,
\be\label{bolt-fourier}
\partial_t \hat f(\xi, t)= \int_{n\in S^2}
\left (\hat f(\xi^+,t) \hat f(\xi^-,t)-\hat f(\xi,t)\hat
f(0,t)\right ) b\left (\frac{\xi}{|\xi|}\cdot n \right ) \, \d n=
\widehat {Q(f,f)} (\xi,t), \ee as was firstly done by Bobylev (see
for instance  \cite{Bo88}). Here we have used the standard
notations
\[
\xi^+ = \frac \xi 2 + \frac {|\xi|}2 n, \quad \xi^- =
\frac \xi 2 - \frac {|\xi|}2 n.
\]
The first results about the non cut-off case for the weak equation
go back to Arkeryd \cite{A81}, in a more general setting. In
\cite{PT}, A.~Pulvirenti and Toscani reformulated the existence
theory starting from the equation in the Fourier variable both
for the cut-off and non cut-off cases. We briefly recall their
method and their result, because they will be useful in the
following. The classical approach is to find a solution of
equation \eqref{bolt-fourier} through a limiting process on the
solutions of a sequence of cut-off approximating problems in the
following way. Let us consider the following sequence of bounded
functions obtained by cutting out the singularity of $b$ at the
origin \be\label{alpha} \bar b_l= \min(b, l), \quad l\in \N, \ee
and let \be\label{beta} b_l^* = \int_{n\in S^2} \bar b_l\left
(\frac{\xi}{|\xi|}\cdot n\right )\, \d n. \ee Then, define
\be\label{gamma}
 \beta_l = \frac {\bar b_l}{b_l^*},
\ee
so that
\[
\int_{n\in S^2} \beta_l\left (\frac{\xi}{|\xi|}\cdot n\right )\, \d n =1,
\]
and then consider the sequence of Cauchy problems
\be\label{bolt-cut-off}
\left\{
\begin{aligned}
&\partial_\tau \varphi_l (\xi,\tau)=
\int_{n\in S^2}
\left (\varphi_l (\xi^+ ,\tau) \varphi_l ( \xi^- ,\tau)-\varphi_l (\xi,\tau)\varphi_l
(0,\tau)\right ) \beta_l\left (\frac{\xi}{|\xi|}\cdot n\right )
 \, \d n,
\\
& \varphi_l (\xi,0)= \hat f_0(\xi).
\end{aligned}
\right .
\ee
A. Pulvirenti and Toscani proved first the existence and uniqueness of a
solution $\varphi_l$ of the Cauchy problems \eqref{bolt-cut-off}.
Then, letting
\[
\hat f_l(\xi, t):= \varphi_l (\xi, b_l^* t),
\]
they proved the convergence in a suitable setting of a subsequence of
$\hat f_l$ to a solution of the Cauchy problem for the non cut-off equation.
More precisely, the result is as follows~:
\begin{theo} [A. Pulvirenti, Toscani \cite{PT}]\label{esist}
We consider an initial datum
$f_0 \geq 0$ satisfying the following assumptions:
\[
\begin{gathered}
\int_{\R^3} f_0(v) \, \d v =1,\quad
\int_{\R^3} f_0(v)\, v_i \,\d v=0,\ i=1,2,3,\quad
\int_{\R^3} f_0(v)|v|^2\, \d v =3,\\
\int_{\R^3} f_0(v) |\log f_0(v)|\, \d v <\infty,
\end{gathered}
\]
and the following Cauchy problem:
\be\label{cauchy-bolt}
\left\{
\begin{aligned}
&\partial_t \hat f(\xi,t)= \int_{n\in S^2}
\left (\hat f( \xi^+,t) \hat f( \xi^- ,t)-\hat f(\xi,t)\hat f(0,t)\right ) b\left (
\frac{\xi}{|\xi|}\cdot n  \right )
 \, \d n, \quad t>0,\\
& \hat f(\xi, 0)= \hat f_0(\xi)
\end{aligned}
\right . \ee where $ b$ is a nonnegative function of
$L^1_{{\rm loc}}([-1,1[)$ satisfying $ b(\cos \theta) = O\left( \frac 1{
(1-\cos \theta)^{\frac 54}}\right),\ \theta \to 0$. Then,
there exists a nonnegative solution $f\in C^1 \left( [0, +\infty),
L^1(\R^3)\right )$ to eq. (\ref{cauchy-bolt}) satisfying for all
$t>0$~:
\[
\begin{gathered}
\int_{\R^3} f(v, t) \, \d v =1,\quad
\int_{\R^3} f(v,t)\, v_i\,\d v=0,\ i=1,2,3, \quad
\int_{\R^3} f(v,t)|v|^2\, \d v =3,\\\
\int_{\R^3} f(v,t) |\log f(v,t)|\, \d v <\infty .
\end{gathered}
\]
Moreover, for all $t>0$, the Fourier transform $\hat f(\cdot, t)$ of the solution
is obtained as the (uniform on compact sets) limit of a
subsequence of the functions
$\varphi_l (\cdot, b_l^* t) \in C^1 \left( [0, +\infty), C_b(\R^3)\right )$,
solutions of the cut-off Cauchy problems \eqref{bolt-cut-off}, which have the
following explicit representation (called Wild's expansion):
\[
\varphi_l (\xi, \tau)= \e^{-\tau} \sum_{k=0}^\infty \varphi_l^{(k)}(\xi)
(1-e^{-\tau})^k,
\]
where
\[
\begin{aligned}
&\varphi_l^{(0)}(\xi) = \hat f_0(\xi),\\
&\varphi_l^{(k+1)}(\xi) = \frac 1{k+1} \sum_{j=0}^k
\int_{n\in S^2}
\varphi_l^{(j)}( \xi^+ ) \varphi_l^{(k-j)} ( \xi^-) \, \beta_l \left
(\frac{\xi}{|\xi|}\cdot n \right )\, \d n.
\end{aligned}
\]
\end{theo}
In \cite{TV}, Toscani and Villani proved that under the hypotheses of
Theorem \ref{esist}, the solution for the non cut-off equation is indeed unique,
whereas the uniqueness of the solution for the cut-off problem
\eqref{bolt-cut-off} was already known.

The question whether an extra property satisfied by the initial
datum $f_0$ propagates along the solution has been already
addressed concerning Sobolev or Lebesgue regularity. In \cite{d-m}
Desvillettes and Mouhot proved the uniform propagation of $L^p$
moments for both the cut-off and non cut-off equations (Cf. also
\cite{gus1}, \cite{gus2}, \cite{mvv} for earlier works on the propagation of $L^p$
regularity).
 In \cite{CGT}, Carlen, Gabetta and Toscani
proved for the cut-off equation that also all the $H^s$ Sobolev
norms remain uniformely bounded if they exist initially (Cf. also \cite{mvv} for related results in the case of
hard potentials and hard spheres). When the non cut-off equation
 is considered,
the same is true, and moreover the $H^s$ norms are immediately created (Cf. \cite{des1, des2}).
 In \cite{U84}, Ukai proved for
both the cut-off and non cut-off equations that a regularity
property of Gevrey type satisfied by the initial datum keeps on
being satisfied at least for a finite time by the solution. We
shall come back to this result later.
\par
Note finally that many papers address the important question of the propagation
of the behavior of the solution with respect to large $v$ (that is, propagation of moments, evolution of Maxwellian tales, etc.).
We do not investigate in this direction in this work.
\bigskip

This paper is devoted to the discussion of the following question: if the initial
datum $f_0$ satisfies the upper bound
\[
|\hat f_0(\xi)| \leq K_1 e^{-K_2|\xi|^s},\quad \xi\in\R^3, \ K_1 \geq 1,\ K_2>0, \ s>0,
\]
does the solution (of the cut-off or non cut-off Boltzmann equation with Maxwellian molecules) keep on satisfying the same property? The answer
is positive, provided that $s \in ]0,2]$ and we allow the constants $K_1, K_2$ to be different from those of
the initial datum.
\par
 More
precisely, the result we are going to prove is the following:
\begin{theo}\label{7}
Let $f_0$ be a nonnegative function satisfying
 \[
\begin{gathered}
\int_{\R^3} f_0(v) \, \d v =1,\quad \int_{\R^3} f_0(v)\, v_i \,\d
v=0,\ i=1,2,3,\quad
\int_{\R^3} f_0(v)|v|^2\, \d v =3,\\
\int_{\R^3} f_0(v) |\log f_0(v)|\, \d v <\infty.
\end{gathered}
\] If $f_0$ is such that
$$
\sup_{{\xi}\in \R^3}|\hat f_0({\xi})|{\rm
e}^{K_2\psi(|{\xi}|^2)}\leq K_1,
$$
for some $K_1\geq 1$, $K_2>0$, and for some concave function
$\psi:[0,+\infty)\to [0,+\infty)$, such that $\psi(0)=0$,
$\psi(r)\leq r$ for $r$ large enough and $\psi(r)\to +\infty $ for
$r\to +\infty$, then there exist $R_0>0$, $K>0$ such that the
unique solution of the Cauchy problem \eqref{cauchy-bolt} with $f_0$
as initial datum satisfies
\[
\begin{aligned}
&\sup_{|{\xi}|< R_0}|\hat f({ \xi},t)|{\rm e}^{K
|{ \xi}|^2}\leq 1,\quad t\geq 0,\\
&\sup_{|{\xi}|\geq R_0}|\hat f({\xi},t)|{\rm e}^{K \psi(|{
\xi}|^2)}\leq 1, \quad t\geq 0.
\end{aligned}
\]
\end{theo}
Denoting by $G^\nu(\R^3)$ the space of Gevrey functions and by
$G^\nu_0(\R^3)$ the space of Gevrey functions with compact support (we
shall recall in Section 4 their definition), we are able to deduce
from the previous result the propagation along the solution of a
Gevrey-type regularity satisfied by the initial datum.
\begin{coro} \label{condiz}
Let  $f_0$ be a nonnegative function satisfying
\be\label{integ}
\begin{gathered}
\int_{\R^3} f_0(v) \, \d v =1,\quad \int_{\R^3} f_0(v) v_i \,\d
v=0,\ i=1, 2, 3,\quad
\int_{\R^3} f_0(v)|v|^2\, \d v =3,\\
\int_{\R^3} f_0(v) |\log f_0(v)|\, \d v <\infty.
\end{gathered}
\ee
\begin{enumerate} [i)]
\item If $\nu>1$ and $f_0\in G^\nu_0(\R^3)$, then the solution
$f(\cdot, t)$ of the Cauchy problem (\ref{cauchy-bolt})
 is in $G^\nu(\R^3)$, uniformly for all $t\geq
0$. \item If $\nu\geq 1$, $f_0\in G^\nu (\R^3) \cap {\cal
S'}(\R^3)$, and moreover satisfies $\sup_{\xi\in \R^3} |\hat
f_0(\xi) | \leq K_1 e^{-K_2 |\xi|^{\frac 1\nu}}$ for $K_1 \geq 1$,
$K_2>0$, then the solution $f(\cdot, t)$ of
the Cauchy problem (\ref{cauchy-bolt})
 is in $G^\nu(\R^3)$,
uniformly for all $t\geq 0$.
\end{enumerate}
\end{coro}
The plan of the paper is the following:
\begin{Myitemize}
\item in Section 2, we shall present the result in a simpler form
and for the so-called Kac model (which is 1-dimensional and
describes radially symmetric solutions of the Boltzmann equation);
\item in Section 3, we shall generalize the result both to
Boltzmann equation and to more general bounds on the initial datum;
\item finally, in Section 4, we shall recall the main definitions
of Gevrey functions, and we shall state the propagation result of
a Gevrey-type regularity.
\end{Myitemize}
{\bf Acknowledgements:} The authors would like to thank G.~Toscani for  useful discussions
about this problem.

\section{The Kac equation}
In this section, we present our result in the simpler
case of the Kac equation. This equation, in its cut-off or
non cut-off version, is
obtained when one considers radially symmetric solutions of the
homogenous Boltzmann equation for Maxwellian molecules. It
reads:
\[
\partial_t f(v,t)= \int_{w\in\R}\int_{\theta \in [-\frac {\pi} 2,
\frac {\pi} 2]} \left (f(\tilde{v}_*, t) f(\tilde{w}_*, t) -f(v,t)f(w,t)\right )
b(\theta) \, \d\theta \, \d w = \tilde{Q}(f,f)(v,t).
\]
Here, $f(v,t) : \R
\times \R^+ \freccia \R$ is the probability density of a gas of
one dimensional particles which depends only on the velocity
$v\in\R$ at the time $t\geq 0$, and which evolves through collisions which
conserve energy but not momentum.
The relations between the velocities $(\tilde{v}_*,\tilde{w}_*)$ of two particles before the
collision and $(v, w)$ after it are the following
\[
\left\{
\begin{aligned}
& \tilde{v}_*= v\cos \theta + w\sin \theta,\\
& \tilde{w}_* = v\sin \theta -w\cos \theta.
\end{aligned}
\right .
\]
We shall make the following non cut-off assumption on the collision kernel $b$:
\[
b(\theta) =  O_{\theta \to 0} \bigg(\frac {\cos \theta}{|\sin \theta|^\gamma} \bigg),\quad \gamma
\in ]1,3[.
\]
Actually, this kind of assumption for the Kac equation was introduced by
Desvillettes in \cite{des1} whereas, in the original equation, $b(\theta)$ is a
strictly positive constant.
In the same way as for the Boltzmann equation, it is useful to consider the
Cauchy problem in the Fourier variable
\be\label{kac-fourier}
\left \{
\begin{aligned}
&\partial_t \hat f(\xi,t)= \int_{\theta \in [-\frac {\pi} 2, \frac {\pi} 2]}
\left (\hat f(\xi \cos \theta, t) \hat f(\xi \sin \theta, t) -\hat
f(\xi,t)\hat f(0,t)\right ) b(\theta) \, \d\theta = \widehat
{\tilde{Q}(f,f)} (\xi,t),\\
&\hat f(\xi, 0)= \hat f_0(\xi),
\end{aligned}
\right .
\ee
where the even initial datum $ f_0 \geq 0$ satisfies the assumptions:
\be\label{in-cond}
\int_{\R} f_0(v) \, \d v =1,
\quad \int_{\R} v^2 f_0(v)\, \d v =1, \quad \int_\R f_0(v)
|\log f_0(v)|\, \d v <\infty.
\ee
The Kac equation shares with the homogenous Boltzmann
equation for Maxwellian molecules the existence and uniqueness
theory for the solutions.
By considering the sequence of cut-off approximating
problems
\be\label{kac-cut-off} \left\{
\begin{aligned}
&\partial_\tau \varphi_l (\xi,\tau)= \int_{\theta \in [-\frac {\pi}
2, \frac {\pi} 2]}
\left (\varphi_l (\xi \cos \theta, \tau) \varphi_l (\xi \sin \theta, \tau) -
\varphi_l (\xi,\tau)\varphi_l (0,\tau)\right ) \beta_l (\theta) \, \d\theta,\\
& \varphi_l (\xi,0)= \hat f_0(\xi),
\end{aligned}
\right .
\ee
where each $\beta_l (\theta)$ is a bounded function defined as in
\eqref{alpha}, \eqref{beta} and \eqref{gamma}, it is possible to prove that
each Cauchy problem \eqref{kac-cut-off} has a unique solution $\varphi_l$, which
has the following explicit representation, called Wild's expansion:
\[
\varphi_l (\xi, \tau)= \e^{-\tau} \sum_{n=0}^\infty
\varphi_l^{(n)}(\xi) (1-e^{-\tau})^n,
\]
where
\[
\begin{aligned}
&\varphi_l^{(0)}(\xi) = \hat f_0(\xi),\\
&\varphi_l^{(n+1)}(\xi) = \frac 1{n+1} \sum_{j=0}^n \int_{\theta \in
[-\frac {\pi} 2, \frac {\pi} 2]}
\varphi_l^{(j)}(\xi\cos \theta) \varphi_l^{(n-j)} (\xi \sin \theta) \,\beta_l
(\theta)\,
\d \theta.
\end{aligned}
\]
Finally, letting
\[
\hat f_l(\xi, t):= \varphi_l (\xi, b_l^* t),
\]
it is possible to establish the (uniform on compact sets) convergence of a
subsequence of $\hat f_l$ to a solution $\hat f$ of
the Cauchy problem for the (non necessarily cut-off) equation
\eqref{kac-fourier}.

\bigskip

Let us suppose now that the initial datum $f_0$
satisfies the extra property:
 \[
|\hat f_0(\xi)|\leq {\e}^{-K|\xi|^s},\quad \xi\in \R,\ K>0, \ s\in (0,2].
\]
Thanks to the representation of
the solution of the cut-off equation \eqref{kac-cut-off} in Wild's expansion,
it is straightforward to prove that the
solution itself satisfies the same upper bound.
Indeed, by a direct computation, we have that $| \varphi_l^{(1)}(\xi)|\leq
{\e}^{-K|\xi|^s}$, since
$$
\varphi_l^{(1)}(\xi){\e}^{K|\xi|^s} = \int_{\theta } {\e}^{K|\xi|^s-K|\xi\cos
\theta|^s-K|\xi\sin \theta|^s}
\hat f_{0}(\xi\cos \theta){\e}^{K|\xi\cos\theta|^s}\hat f_0(\xi
\sin\theta){\e}^{K|\xi\sin\theta|^s}\,\beta_l(\theta)\, \d \theta,
$$
and $1-|\cos \theta|^s-|\sin \theta|^s\leq 0$, for $s\in (0,2]$.
Hence by an immediate iteration argument, the same inequality
holds for any $\varphi_l^{(n)}(\xi)$, and finally for the solution
of the cut-off equation for any $t\geq 0$. Passing to the limit
when $l\to +\infty$ in the estimate $\varphi_l (\xi, b_l^* t) \le
{\e}^{-K|\xi|^s} $, we see that the inequality also holds for the
solution of the (non necessarily cut-off) equation \eqref{kac-fourier}.

Due to the non-linearity of the collision operator, if we now
consider the weaker assumption
\[
|\hat f_0(\xi)|\leq
 K_1{\e}^{-K_2|\xi|^s}, \quad \xi\in \R,\, K_1>1,\, K_2>0,\, s\in (0,2],
\]
the same argument allows to prove that
the solution of each cut-off equation \eqref{kac-cut-off} satisfies the same
upper bound, but only for a finite interval of time. In this case, by letting
$l$ go to infinity, the interval of time where the
estimate is true can reduce to nothing.

In spite of this, we prove in this section that the condition
$|\hat f_0(\xi)|\leq K_1{\e}^{-K_2|\xi|^s}$ propagates
(though possibly with different constants $K_1$ and $K_2$ )
along the solution of the (non necessarily cut-off) Kac equation.

\subsection{Some preliminary properties of initial data}
In this section we emphasize some useful properties satisfied by
any even, nonnegative function $g$ such that $\int_{\R}g(v)\, \d
v=1$ and $\int_{\R}v^2\, g(v)\, \d v=1$.
\begin{lem}\label{lemma3-c}
Let $g$ be a nonnegative, even function, satisfying
\[
\int_{\R} g(v) \, \d v =1,\quad \int_{\R} g(v)\, v^2 \, \d v
=1.
\]
Then, there exist $\rho >0$ and $\tilde K>0$ such that for all
$|\xi|\leq \rho$:
\[
|\hat g(\xi)| \leq \e^{-\tilde K |\xi|^2}.
\]

\end{lem}

\begin{proof}
We observe that under the hypotheses of the lemma, $\hat{g}$ is of
class $C^2$ and satisfies the following property
: $\hat{g}(0) =
1$, $\hat{g}'(0) = 0$, and $\hat{g}''(0) = -1$. Using a Taylor
expansion of $\hat{g}$ at order $2$, we obtain $\hat{g}(\xi) = 1 -
\frac 12 \xi^2 + o(\xi^2)$ when $\xi \to 0$. Then, the estimate of
the lemma holds for any $\tilde K \in ]0, \frac 12[$.
\end{proof}

Then, we prove the:
\begin{lem}\label{lemma1-c}
Let $g \geq 0$ such that $\int_\R g(v) \, \d v =1$. Then, for all $r>0$, there
exist $C_r\in(0,\frac 12)$ and $\tilde C_r\in(0,\frac 12)$ such that
\[
\begin{aligned}
&\int_\R g(v) \sin^2\left(\frac {v\xi}2\right ) \, \d v \geq C_r,\quad |\xi| >r,\\
&\int_\R g(v) \cos^2\left(\frac {v\xi}2\right ) \, \d v \geq \tilde C_r, \quad |\xi|
>r.
\end{aligned}
\]
\end{lem}
\begin{proof}
We only prove the first inequality, since the second one can be proven in
exactly the same way.

Thanks to Lebesgue's dominated convergence theorem and thanks to the absolute continuity of the
measure
$\nu(E):= \int_E g(v)\, \d v$ with respect to the Lebesgue measure, there exist $R>0$
and $\delta>0$ such that
for all measurable set $A \subset \R$ such that $|A| \le \delta$, we have
\be\label{prima}
\int_{A^c \cap B(0,R)} g(v)\, \d v \geq \frac 12.
\ee
Let $\xi \in \R$ be fixed. For $\mu \in ]0, \pi/2[$, we define
\[
K_{\mu, R} := \left \{v\in \R, |v| \leq R\ {\rm and\ }\exists\ k\in\Z, \left| \frac
{v\xi}2 -k\pi\right | \leq \mu \right \}.
\]
It is clear that
\[
\left| K_{\mu, R}\right | \leq \left( \frac {|\xi|R}\pi +1\right ) \frac {4\mu}
{|\xi|} = 4\mu\left ( \frac R \pi +\frac 1{|\xi|}\right )
\]
so that, when $|\xi| \geq r$ we have $ \left| K_{\mu, R}\right | \leq 4\mu\left (
\frac R \pi +\frac 1r\right )$.
When $\mu= \frac {\delta}{4\left( \frac R\pi +\frac 1r\right)}$, we see
thanks to \eqref{prima} that
\[
\int_{K_{\mu,R}^c\cap B(0,R)} g(v)\, \d v \geq \frac 12.
\]
We can therefore conclude that
\[
\int_\R g(v) \sin^2\left(\frac {v\xi}2\right ) \, \d v \geq \int_{K_{\mu,R}^c\cap
B(0,R)} g(v)\sin^2\left(\frac {v\xi}2\right ) \, \d v \geq \frac 12 \sin^2 \mu :=C_r.
\]
\end{proof}

\subsection{The propagation theorem}
We are now in position to state the theorem.
\begin{theo}\label{5}
Let $f_0$ be a nonnegative, even function, satisfying
\[
\int_{\R} f_0(v) \, \d v =1,
\quad
\int_{\R} f_0(v)\, v^2\, \d v =1,\quad
\int_{\R} f_0(v) |\log f_0(v)|\, \d v <\infty.
\]
We suppose that $f_0$ is such that
$$
\sup_{\xi\in \R}|\hat f_0(\xi)|{\rm e}^{K_2|\xi|^s}\leq K_1
$$
for some $K_1\geq 1$, $K_2>0$, and $0<s\leq 2$. Then there exist
$R_0>0$, $K>0$ such that the unique solution of the Cauchy problem
\eqref{kac-fourier}
satisfies:
\be\label{finale}
\begin{aligned}
&\sup_{|\xi|< R_0}|\hat f(\xi,t)|{\rm e}^{K
|\xi|^2}\leq 1,\quad t\geq 0,\\
&\sup_{|\xi|\geq R_0}|\hat f(\xi,t)|{\rm e}^{K
|\xi|^s}\leq 1, \quad t\geq 0.
\end{aligned}
\ee
\end{theo}
We begin by proving the following proposition.
\begin{prop}\label{espo}
Let $g$ be a nonnegative, even function, satisfying
\[
\int_\R g(v) \, \d v =1,\quad \int_\R g(v)\, v^2 \, \d v =1.
\]
Let us suppose moreover that, for given $s\in (0, 2]$, $K_1 >
1$ and $K_2>0$, $g$ satisfies the following bound:
\[
|\hat g(\xi)| \leq K_1\e^{-K_2 |\xi|^s}, \quad \xi\in \R.
\]
Then, there exists $\eta >0$ such that for all $R > \eta $, there
exists $K>0$ (depending on $R$) such that
\[
|\hat g(\xi)| \leq \begin{cases}
\e^{-K |\xi|^2}, & |\xi| < R,\\

\e^{-K |\xi|^s}, & |\xi| \geq R.
\end{cases}
\]
\end{prop}

\begin{proof}
We have already proven in Lemma \ref{lemma3-c} that $|\hat g(\xi)|
\leq \e^{- \tilde K |\xi|^2}$ for $|\xi| \leq \rho$, where
$\tilde K$ and $\rho$ are suitably chosen. Now, let
$\eta=\left(\frac{\log {K_1}}{K_2}\right)^{\frac 1s}$. For every
$R>\eta$, we can find $0<K_3 <K_2$ such that
\[
K_1\e^{-K_2 |\xi|^s} \leq \e^{-K_3 |\xi|^s}, \quad |\xi| \geq R,
\]
so that
\[
|\hat g(\xi)| \leq \e^{-K_3 |\xi|^s}, \quad |\xi| \geq R.
\]
It is now enough to find $K_4 >0$ such that
\[
|\hat g(\xi)| \leq \e^{-K_4 |\xi|^2}, \quad \rho <|\xi| < R.
\]
Since $g$ is an even function, we have
\[
\hat g(\xi) = \int_{\R} g(v) \e^{-i \xi v} \, \d v= \int_\R g(v)
\left( \frac{\e^{-i\xi v} +\e^{i\xi v}} 2\right )\, \d v =
\int_{\R} g(v) \cos ( \xi v) \, \d v.
\]
Then, $\hat g$ is real and $|\hat g (\xi)| \leq 1$ for all $\xi
\in \R$. Moreover
\[
\begin{aligned}
&1-\hat g(\xi) = \int_\R g(v) (1-\cos (\xi v) )\, \d v = 2\int_\R g(v) \sin^2
\left(\frac {\xi v} 2\right )\, \d v,\\
&\hat g(\xi) +1 = \int_\R g(v) (1+\cos (\xi v) )\, \d v = 2\int_\R
g(v) \cos^2 \left(\frac {\xi v} 2\right )\, \d v.
\end{aligned}
\]
According to Lemma \ref{lemma1-c}, we know that
\[
\begin{aligned}
&1-\hat g(\xi) \geq 2\,C_\rho,\quad |\xi|> \rho,\\
&\hat g(\xi) +1\geq 2\,\tilde{C_\rho},\quad |\xi| >\rho.
\end{aligned}
\]
Therefore,
\[
|\hat g(\xi)| \leq 1-\min(2\, C_\rho,2\, \tilde C_\rho),\quad |\xi| > \rho.
\]
Now, there exists $K_4>0$ such that
\[
1-\min(2\,C_\rho,2\,\tilde{C_\rho}) \leq \e^{-K_4 R^2},
\]
which implies
\[
 |\hat g(\xi)| \leq 1-\min(2\,C_\rho,2\,\tilde{C_\rho}) \leq \e^{-K_4 R^2} \leq
\e^{-K_4 |\xi|^2}, \quad \rho < |\xi| \leq R.
\]
We can conclude letting $K= \min(\tilde K, K_3, K_4)$.
\end{proof}

\noindent {\bf Proof of Theorem \ref{5}: the cut off case.}
Thanks to Proposition \ref{espo}, there exists $\eta >0$ such that for any
$R_0>\eta$, there exists a strictly positive $K$ such that the initial datum
$\hat f_0$ satisfies
\begin{equation}\label{(tre)}
\begin{aligned}
&\sup_{|\xi|< R_0}|\hat f_0(\xi)|{\rm e}^{K
|\xi|^2}\leq 1,\\
&\sup_{|\xi|\geq R_0}|\hat f_0(\xi)|{\rm e}^{K
|\xi|^s}\leq 1.\\
\end{aligned}
\end{equation}
In order to prove the theorem for the cut-off case, it is enough to
establish that any $\varphi_l^{(n)} $ in Wild's sums satisfies
\eqref{(tre)}. Let us check that this is true for $\varphi_l^{(1)}$.
Let us define
\[
H(|\xi|)=\begin{cases}
{K |\xi|^2}, &|\xi|< R_0,\\
{K |\xi|^ s}, & |\xi|\geq R_0.
\end{cases}
\]
Condition \eqref{(tre)} on the initial datum $f_0$ reads therefore
\[
\sup_{\xi \in \R}|\hat f_0(\xi)|{\rm e}^{H(|\xi|)}\leq 1.
\]
Then
\[
\begin{aligned}
& \left| {\rm e}^{H(|\xi|)}\varphi_l^{(1)}(\xi)\right|\\
\leq & \int_{\theta \in [-\frac \pi 2, \frac \pi 2]}{\rm
e}^{H(|\xi|)-H(|\xi\sin \theta|)-H(|\xi\cos \theta|)} {\rm
e}^{H(|\xi\sin \theta|)} |\hat f_0(\xi\sin \theta)|{\rm
e}^{H(|\xi\cos \theta|) }|\hat f_0(\xi\cos \theta)|
\beta_l(\theta)\d \theta\\
\leq &\int_{\theta \in [-\frac \pi 2, \frac \pi 2]} {\rm
e}^{H(|\xi|)-H(|\xi\sin \theta|)-H(|\xi\cos \theta|)}
\beta_l(\theta)\d \theta.
\end{aligned}
\]
Since $\int_\theta \beta_l(\theta) \d \theta=1$, we end the estimate
by proving that $H(|\xi|)-H(|\xi\sin \theta|)-H(|\xi\cos
\theta|)\leq 0$ for $\xi\in\R$ and $\theta \in [-\frac {\pi}2,
\frac {\pi}2]$ if $R_0 \ge 1$. Thanks to the symmetries
of the function $H$ with respect to $\theta$, we can restrict
ourselves to the interval $[0,\frac \pi 4]$. Now, when $|\xi|<
R_0$, we have
$$
H(|\xi|)-H(|\xi\sin \theta|)-H(|\xi\cos \theta|)=K |\xi|^2(1-(\sin
\theta)^2-(\cos \theta)^2)=0.
$$
If $|\xi|\geq R_0$ and $|\xi\sin \theta|\geq R_0$, $|\xi\cos
\theta|\geq R_0$, then
$$
H(|\xi|)-H(|\xi\sin \theta|)-H(|\xi\cos \theta|)=K |\xi|^s(1-(\sin
\theta)^s-(\cos \theta)^s)\leq 0
$$
for $0<s\leq 2$. Whenever $|\xi|\geq R_0$ and $|\xi\sin \theta|<
R_0$, $|\xi\cos \theta|< R_0$ we have
$$
H(|\xi|)-H(|\xi\sin \theta|)-H(|\xi\cos
\theta|)=K \left(|\xi|^s-|\xi|^2\left((\sin \theta)^2+(\cos
\theta)^2\right)\right )=K \left(|\xi|^s-|\xi|^2\right ).
$$
If we choose $R_0\geq 1$, we can conclude since $|\xi|\geq R_0$ that
\[
|\xi|^s-|\xi|^2\leq 0.
\]
If now
$|\xi|\geq R_0$ and $|\xi\sin \theta|< R_0$, $|\xi\cos \theta|\geq
R_0$, we have
$$
H(|\xi|)-H(|\xi\sin \theta|)-H(|\xi\cos
\theta|)=K \left(|\xi|^s-|\xi|^2(\sin \theta)^2-|\xi|^s(\cos \theta)^s\right)
$$
$$ \le K\, \left( |\xi|^s - |\xi|^s\,(\sin\theta)^2 - |\xi|^s\,(\cos\theta)^2 \right)=0.
$$
Note that since $0\le \theta \le \pi/4$, there is no other case to treat.

\vskip 0.5truecm \noindent {\bf The non cut-off case.} As we have
recalled in the introduction of section 2, the solution of \eqref{kac-fourier}
in the non cut-off case is obtained as the limit of a subsequence
 of the solutions of the Cauchy problems \eqref{kac-cut-off}. Since the
estimate on $\varphi_l (\xi, \tau)$ holds true (for any $\tau\geq 0$), the
same is valid for $\hat f_l(\xi,t)$ and hence for $\hat f(\xi,t)$.
\hfill $\scriptstyle\square$

\rem In Theorem \ref{5}, the hypothesis that $f_0$ is even could be
replaced by the weaker hypothesis that $\int_\R f_0(v)\,v\, \d v=0$. Since the Kac
equation comes from the Boltzmann equation when one considers
radially symmetric solutions, it is however natural to study only even
initial data.

\section{Propagation for the Boltzmann equation}
We would like to extend to the solution of the Boltzmann equation
\eqref{bolt-fourier} the results proven in the previous
section for the solution of the Kac equation. Two kinds of
extensions are in order: first, we have to pass from the
one-dimensional to the three-dimensional setting; second we
would like to state the result considering not only functions like
$\e ^ {-|\xi|^s}$, but also like $\e^{-\psi(|\xi|^2)}$, where $\psi$
is a suitable concave function. We now begin by restating the
lemmas of the previous section in three dimensions. We shall only indicate the major
modifications in the proofs.

\begin{lem}\label{lemma3-b}
Let $g : \R^3 \to \R$ be a nonnegative function satisfying
\[
\int_{\R^3} g(v) \, \d v =1,\quad \int_{\R^3} g(v) \,v_i \d v =0, \ i=1,2,3, \quad
\int_{\R^3} g(v)\, |v|^2 \, \d v =3.
\]
Then
\, there exist $\rho >0$ and $\tilde K>0$ such that for all $|\xi|\leq \rho$:
\[
|\hat g(\xi)| \leq \e^{-\tilde K |\xi|^2}.
\]
\end{lem}

\noindent
\begin{proof}
We observe that the result of the lemma is not changed when $g$ is
replaced by $g \circ R$, where $R$ is any rotation of $\R^3$. As a
consequence, we can suppose that the symmetric matrix
$(\int_{\R^3} g(v)\,v_i\,v_j\, \d v)_{i,j \in \{1,2,3\}} $ is
diagonal. Moreover, since $g \in L^1(\R^3)$ and $ \int_{v\in \R^3}
g(v)\,|v|^2\, \d v = 3$, we see that $\int_{v\in \R^3} g(v)\,
v_i^2\, \d v >0$ for $i=1,2,3$. \par Then, $\hat g(\xi) = 1 -
\sum_{j=1}^3 \lambda_j \,\xi_j^2 + o_{\xi \to 0}(|\xi|^2)$, with
$\lambda_j > 0$ for $j=1,2,3$, and we conclude like in Lemma
\ref{lemma3-c}.
\end{proof}

\begin{lem}\label{lemma1-b}
Let $g : \R^3 \to \R$ be a nonnegative function satisfying
$\int_{\R^3} g(v) \, \d v =1$. Then, for all $r>0$, there exists
$C_r\in(0,\frac 12)$ such that for all $\theta \in\R$,
$$
\int_{\R^3} g(v) \sin^2\left(\frac {v\cdot \xi + \theta}2\right ) \, \d v
 \geq C_r,\quad |\xi| >r .
$$
\end{lem}
\begin{proof}
The proof follows the same lines as that of Lemma \ref{lemma1-c}. Let $\xi \in \R^3$ be fixed. We start by choosing in
$\R^3$ an orthogonal system in which the unitary vector along the z-axis is
$\frac \xi {|\xi|}$. As in the one-dimensional case, there exist
 $R>0$ and $\delta>0$ such that
for all measurable set $A \subset \R^3$ such that $|A| < \delta$,
we have
\[
 \int_{A^c \cap Q(0,R)} g(v)\, \d v \geq
\frac 12,
\]
where $Q(0,R)$ is the cube centered at the origin:
\[
Q(0,R)= \{v\in \R^3:\, |v_i| \leq R,\, i=1,2,3\}.
\]
Now, for $\mu \in\ ]0, \frac \pi 2[$, we define
\[
K_{\mu, R, \theta} := \left \{v\in \R^3, v\in Q(0,R)\ {\rm and\ }\exists\
k\in\Z, \left| \frac {v\cdot \xi + \theta}2 -k\pi\right | \leq \mu \right
\}.
\]
Thanks to the choice of the coordinate system, we have
\[
K_{\mu, R, \theta} = \left\{v\in \R^3, |v_i| \leq R, i=1,2,3\ {\rm and\ }
\exists\ k\in\Z, \left| \frac {v_3|\xi| + \theta}2 -k\pi\right | \leq \mu
\right \}.
\]
So, it is easy to see that
\[
\left| K_{\mu, R, \theta}\right | \leq \left( \frac {|\xi|R}\pi +1\right )
\frac {4\mu} {|\xi|} R^{2} = 4\mu R^{2}\left ( \frac R \pi +\frac
1{|\xi|}\right ),
\]
and we can conclude as in the one-dimensional case.
\end{proof}

\medskip\noindent
We are now in position to prove the main theorem of our paper,
namely Theorem \ref{7}. The proof of this theorem relies on the
following proposition:
\begin{prop}\label{espo-b}
Let $g$ be a nonnegative function satisfying
\[
\int_{\R^3} g(v) \, \d v =1,\quad \int_{\R^3} g(v) v_i \,\d
v=0,\ i=1,2,3, \quad \int_{\R^3} g(v)\, |v|^2 \, \d v =3.
\]
Let us suppose moreover that
\[
|\hat g(\xi)| \leq K_1\e^{- K_2 \phi( |\xi|)}, \quad \xi\in \R^3,
\]
where $K_1\geq 1$, $K_2>0$ and $\phi:[0,+\infty)\to [0,+\infty)$
satisfies $\lim_{t\to +\infty}
\phi(t) = +\infty$. Then, there exists $\eta >0$ such that for all
$R > \eta $, there exists $ K>0$ (depending on $R$) such that
\[
|\hat g(\xi)| \leq \begin{cases}
\e^{- K |\xi|^2}, & |\xi| < R,\\
\e^{- K \phi(|\xi|)}, & |\xi| \geq R.
\end{cases}
\]
\end{prop}
\begin{proof}
The proof of this proposition
is only a slight modification of that of Proposition \ref{espo}.
We only point out the few differences.
\par
First, we explain how to find $\eta >0$: we can fix $K_3 \in (0,
K_2)$ and let $\eta = \eta (K_3)$ be a positive constant such that
$\varphi(|\xi|)\geq \frac{1}{K_2-K_3}\log K_1$ for every
$|\xi|\geq \eta$; then of course for every $R\geq\eta$,
\[
K_1\e^{-K_2 \varphi(|\xi|)} \leq \e^{-K_3 \varphi(|\xi|)}, \quad
|\xi| \geq R.
\]
\par
Second, we observe that for all $\xi \in \R^3$ it is possible to find $\theta\in\R$
(depending on $\xi$)  such that
$$  |\hat{g}(\xi)| =  \hat{g}(\xi)\,e^{i\theta} = \int_{\R^3} g(v)\, \cos \left(
\xi\cdot v + \theta\right ) \, \d v.$$
So,  we have
$$
1-|\hat{g}(\xi)| =  2\int_{\R^3} g(v) \sin^2\left(\frac{\xi\cdot v + \theta}2\right)
\, \d v .
$$
Thanks to Lemma \ref{lemma1-b}, there exists $C_\rho$ such that
$$
2\int_{\R^3} g(v)\,
\sin^2\left(\frac{\xi\cdot v + \theta}2\right) \, \d v \ge 2\,C_\rho,\quad |\xi| >
\rho.
 $$
Then, we can conclude as in
the proof of Proposition \ref{espo}.
\end{proof}

\noindent {\bf Proof of Theorem \ref{7}}.
As we did for Kac equation, for each of the Cauchy problems
\eqref{bolt-cut-off}, we write the solution $\phi_l$ under the form of a Wild's
expansion.
\par
In order to prove the bound for the solution,
it is enough to prove it for every term $\varphi_l^{(n)}$ in the sum.
We define
\[
H(|\xi|)=\begin{cases}
{K |\xi|^2}, &|\xi|< R_0,\\
{K \psi(|\xi|^ 2)}, & |\xi|\geq R_0,
\end{cases}
\]
where $R_0$ will be chosen (large enough) later, and $K$ is given
by Proposition \ref{espo-b}. Thanks to this proposition, the
initial datum satisfies
\[
\sup_{\xi \in \R^3}|\hat f_0(\xi)|{\rm e}^{H(|\xi|)}\leq 1.
\]
We recall the identities (\cite{des2}, page 56)
\[
\begin{aligned}
&\left|\xi^+\right |= |\xi| \cos \frac \theta 2,\\
&\left|\xi^-\right | = |\xi| \sin \frac \theta 2.
\end{aligned}
\]
For the first term $\varphi_l^{(1)}$ we have:
\[
\begin{aligned}
\left|\varphi_l^{(1)}({\xi}){\rm e}^{\;H(|{\xi}|)}\right| &=
\left|\int_{S^2} {\rm e}^{\;H(|{\xi}|)
-\;H ( |\xi^+|)
-\; H(|\xi^-|) }
\hat f_0(\xi^+)
{\rm e}^{\;H(|\xi^+ | )}
\hat f_0(\xi^-)
{\rm e}^{\;H(|\xi^-| )}
 \beta_l \left( \frac \xi {|\xi|} \cdot n\right )
  \, \d n\right |\\
&\leq \int_{\theta\in (0, \pi )}\int_{\varphi \in (0,2\pi)}{\rm
e}^{\;H(|{ \xi}|)- \; H(|{\xi}| \cos\frac{\theta}{2})- \;
H(|{\xi}| \sin\frac{\theta}{2})} \beta_l(\cos \theta) \sin
\theta\, \d\theta\, \d\varphi,
\end{aligned}
\]
where in the last integral, we have used the spherical coordinates with
$\frac{\xi}{|\xi|}$ as $z$-axis. Then, in order to establish
$\left|\varphi^{(1)}_l({\xi}){\rm e}^{\;H(|{ \xi}|)}\right|\leq 1$,
we show that for $R_0$
large enough,
\[
\;H(|{\xi}|)- \; H\left(\left|{\xi}
\cos\frac{\theta}{2}\right|\right)- \;
H\left(\left|{\xi}\sin\frac{\theta}{2}\right|\right) \leq 0,\quad \theta \in (0,
{\pi}),\ \xi \in \R^3.
\]
We denote $\tilde \theta=\frac{\theta}{2}$.
Thanks to the symmetries of the functions $H(|\xi|)-H(|\xi\sin \tilde\theta|)-
H(|\xi\cos \tilde\theta|)$
with respect to $\tilde\theta$, we can restrict ourselves to the interval
$(0,\frac \pi 4)$. Now, the case $|\xi|< R_0$ is the same as in
Theorem \ref{5}. If $|\xi|\geq R_0$ and $|\xi\sin
\tilde\theta|\geq R_0$, $|\xi\cos \tilde\theta|\geq R_0$, then
$$
H(|\xi|)-H(|\xi\sin \tilde\theta|)-H(|\xi\cos \tilde\theta|)=K
(\psi(|\xi|^2)-\psi(|\xi\sin \tilde\theta|^2)-\psi(|\xi\cos
\tilde\theta|^2)).
$$
Thanks to the concavity property of $\psi$ and the fact that
$\psi(0)=0$, we have $\psi(|\xi|^2(\sin\tilde\theta)^2)\geq
(\sin\tilde\theta)^2\psi(|\xi|^2)$ and
$\psi(|\xi|^2(\cos\tilde\theta)^2)\geq
(\cos\tilde\theta)^2\psi(|\xi|^2)$. Hence
$$
\psi(|\xi|^2)-\psi(|\xi\sin\tilde {\theta}|^2)-\psi(|\xi\cos
\tilde\theta|^2)\leq \psi(|\xi|^2)-(\sin
\tilde\theta)^2\psi(|\xi|^2)-(\cos\tilde{ \theta})^2\psi(|\xi|^2)=
0.
$$
Whenever $|\xi|\geq R_0$ and
$|\xi\sin\tilde\theta|< R_0$, $|\xi\cos\tilde\theta|< R_0$ we have
$$
H(|\xi|)-H(|\xi\sin\tilde \theta|)-H(|\xi\cos\tilde\theta |)=K
(\psi(|\xi|^2)-|\xi|^2((\sin \tilde\theta)^2+(\cos
\tilde\theta)^2))=K (\psi(|\xi|^2)-|\xi|^2).
$$
If we choose $R_0$ large enough, thanks to the assumption that
$\psi(r) \le r$ for $r$ large enough, and since $|\xi|\geq R_0$,
we can conclude
\[
\psi(|\xi|^2)-|\xi|^2\leq 0.
\]
By using now the concavity property of $\psi$ and the
fact that
$\psi(r)\leq r$, if now $|\xi|\geq R_0$ and
$|\xi\sin\tilde\theta |< R_0$, $|\xi\cos \tilde\theta|\geq R_0$ we
have
$$
\begin{aligned}
H(|\xi|)-H(|\xi\sin \tilde\theta|)-H(|\xi\cos \tilde\theta|)&=K
(\psi(|\xi|^2)-|\xi|^2(\sin\tilde\theta
)^2-\psi(|\xi\cos\tilde\theta
|^2))\\
&\leq K (\psi(|\xi|^2)  - \psi(|\xi|^2)\,(\sin \tilde\theta)^2 -
\psi(|\xi\cos\tilde\theta
|^2))\\
&\leq K\psi(|\xi|^2)(1-(\sin
\tilde\theta)^2-(\cos\tilde\theta)^2)=0.
\end{aligned}
$$
We end the proof by first noticing that a simple induction shows the estimate
$\left|\varphi^{(n)}_l({\xi}){\rm e}^{\;H(|{ \xi}|)}\right|\leq 1$ when $n\ge 1$, and then we may pass to the limit
when $l\to +\infty$ if necessary (that is, in the non cut-off case).
\hfill $\scriptstyle\square$

\rem We would like to point out that the assumption of
concavity for the function $\psi(|\xi|^2)$ is not mandatory and
that the argument exploited in the proof of Theorem \ref{7} could
work also in a more general framework.
\par
By analysing the proof, one can see that, instead of assuming that
$\psi$ is concave, it is in fact enough to assume that for some
$R_0$ large enough,
$$ \psi(\lambda^2\, |\xi|^2) \ge \lambda^2\,\psi(|\xi|^2)$$
when $0 \le \lambda \le 1$, $\lambda\,|\xi| \ge R_0$.
This is true for example when $\psi(t) = \frac12\, \sqrt{t}\, |\log t|$.

\section{Gevrey spaces}
In this section, we translate the propagation result obtained in the previous
sections in terms of Gevrey regularity for the solutions of 
Boltzmann equation. Let us begin by recalling the classical definitions of
Gevrey functions and a useful characterisation of these functions through their Fourier transform.
For more information, the interested reader can consult for instance the book
by Rodino \cite{rod}, from where we have taken the following recalls.
\bigskip

Let $\Omega \subseteq \R^n$ be an open set and let $\nu\geq 1$ be a fixed real
number.
\begin{defi} \label{defgevrey-b}
The class $G^\nu(\Omega)$ of Gevrey functions of order $\nu$ in $\Omega$ is the set of
functions $f \in C^\infty (\Omega)$ satisfying the following property:
for every compact subset $K$ of $\Omega$,
there exists a positive constant $C=C(K)$ such that for all $l \in \N^n$ and
all $x\in K$,
\be \label{defgevrey}
|\partial^l f(x)| \leq C^{|l|+1} (l !)^\nu .
\ee
\end{defi}
Assumption \eqref{defgevrey} can be replaced by other equivalent assumptions,
for example
\[
|\partial^l f(x)| \leq R C^{|l|} (l !)^\nu,
\]
where $R$ and $C$ are two positive constants independent of $l$ and $x\in K$.
It is easy to recognize that $G^1(\Omega)= A(\Omega)$,
 the space of all analytic functions in $\Omega$, and that for
$\nu\leq \tau$, one has $G^\nu(\Omega) \subseteq G^\tau(\Omega)$.
Moreover, it is interesting to underline the following inclusions:
\[
A(\Omega) \subset \bigcap_{\nu>1} G^\nu(\Omega),\quad \bigcup_{\nu\geq 1}
G^\nu(\Omega) \subset C^\infty(\Omega),
\]
which are strict in both cases.
We also recall that the Gevrey class $G^\nu(\Omega)$ is closed under
diffe\-rentiation.

In what follows, we shall also need the following
\begin{defi}\label{gev-comp}
Assume $\nu>1$. We shall denote by $G^\nu_0 (\Omega)$ the vector space of all
$f\in G^\nu(\Omega)$ with compact support in $\Omega$.
\end{defi}
The exclusion of $\nu=1$ in the previous definition is mandatory,
because there are no analytic test functions other than the zero
function. As for the other values $\nu>1$, one could wonder
whether such compact supported functions do exist. An example in
$\R$ is the following: let $r>0$, $\nu >1$, $d= \frac 1 {1-\nu}$
and
\[
\varphi(t)=
\begin{cases}
\e^{-t^d} & t>0,\\
0 & t\leq 0.
\end{cases}
\]
The function
\[
f(x)= \varphi (x+r) \varphi (x-r)
\]
is then in $G^\nu_0(\R)$ (\cite {rod}).

The result that we are going to use in order to relate our propagation result
to Gevrey regularity is the following.
\begin{theo}[\cite{rod}, Theorem 1.6.1 page 31]\label{caratt}
\begin{enumerate}[i)]
\item Let $\nu>1$. If $\phi \in G^\nu_0(\R^n)$, then there exist
positive constants $C$ and $\eps$ such that \be\label{diret} |\hat
\phi (\xi)| \leq C \e^ {-\eps|\xi|^{\frac 1\nu}}, \quad
\xi\in\R^n. \ee \item Let $\nu\geq 1$. If the Fourier transform of
$\phi\in\S'(\R^n)$ satisfies \eqref{diret}, then $\phi\in
G^\nu(\R^n)$.
\end{enumerate}
\end{theo}
We can therefore deduce Corollary \ref{condiz} concerning the
regularity of the solutions of Boltzmann equation.

\noindent {\bf Proof of Corollary \ref{condiz}}. The result is
straightforward from Theorem \ref{caratt} and Theorem \ref{7}. It
is enough to notice that one can replace the uniform estimate
\[
\begin{aligned}
&\sup_{|\xi|< R_0}|\hat f(\xi,t)|{\rm e}^{K
|\xi|^2}\leq 1,\quad t\geq 0,\\
&\sup_{|\xi|\geq R_0}|\hat f(\xi,t)|{\rm e}^{K
|\xi|^s}\leq 1, \quad t\geq 0,
\end{aligned}
\]
obtained in Theorem \ref{7}, by the following uniform estimate:
\[
\sup_{\xi\in \R^3}|\hat f(\xi,t)|{\rm e}^{\tilde K_2 |\xi|^s}\leq
\tilde K_1,\quad t\geq 0,
\]
for $\tilde K_1\geq 1$ and $\tilde K_2 >0$ properly chosen.
Letting now $\nu=\frac 1 s$, one immediately gets the result. This
ends the proof.\hfill $\scriptstyle\square$

\bigskip
 We end up this section
by comparing the regularity result we have just obtained with the
one obtained by Ukai in \cite{U84}. In his work, he considered
among others a Cauchy problem for the homogeneous Boltzmann
equation for Maxwellian molecules both in the cut-off and non
cut-off settings. He considered only initial data $f_0$ satisfying
a strong regularity assumption: for $\alpha \geq 0$, $\rho\geq 0$
and $\nu \geq 1$ he supposed $f_0\in \gamma_{\alpha, \rho}^\nu$,
where
\[
\gamma_{\alpha, \rho}^\nu = \{ g: \|g\|_{\alpha, \rho, \nu}=
\sum_{l\in \N^3} \frac{\rho^{|l|}}{(l!)^\nu} \sup_{v\in\R^3} \e
^{\alpha(1+|v|^2)^{\frac 12}} |\partial^l g(v)| <\infty\}.
\]
Comparing this space with the spaces in  Definition \ref{defgevrey-b}, one
can deduce that initial data in Ukai setting are indeed in the
Gevrey space $G^\nu (\R^3)$, but also decay very strongly at
infinity together with all their derivatives. Ukai was able to
prove by a fixed point argument that there exists a unique, local
in time solution $f$ belonging at every time to a functional space
of the same kind as the initial datum but in which the indices
change with $t$. More precisely, he proved that there exist $T>0$,
$\beta >0$, $\sigma>0$ such that \be\label{ukai-2} \|f(\cdot,
t)\|_{\alpha -t \beta , \rho -t \sigma , \nu} \leq 2\, \|
f_0\|_{\alpha, \rho, \nu}, \quad t \in [0, T]. \ee Since initial
data which belong to $G^\nu_0(\R^3)$ also belong to Ukai space,
the question of comparing the two results is meaningful. Let us
consider a nonnegative  function $f_0 \in G^\nu_0(\R^3)$. If we
suppose moreover that $f_0$ satisfies assumptions \eqref{integ},
then by Theorem \ref{esist} we know that there is a solution $f\in
C^1 \left( [0, +\infty), L^1(\R^3)\right )$ which is unique by
Toscani--Villani's result. Since $G^\nu_0(\R^3) \subset
\gamma_{\alpha, \rho}^\nu$ for all $\alpha\geq 0$, we can deduce
from Ukai result that for all $\alpha \geq 0$ there is a time
$T=T(\alpha)>0$ (possibily finite) such that this solution stays
in the class \eqref{ukai-2} for $t\in [0,T(\alpha)]$ but this
space is not uniform in time (in addition to the fact that for all
$\alpha$ it is difficult to compare the life times $T(\alpha)$).
Our result says instead that the solution stays for $t\in
[0,\infty)$ in the same Gevrey class as its initial datum, without
any information about the decay at infinity and moreover that all
the estimates on the Gevrey seminorms are uniform in time.


\begin{thebibliography}{CGT99}

\bibitem[Ark81]{A81}
L.~Arkeryd.
\newblock Intermolecular forces of infinite range and the {B}oltzmann equation.
\newblock {\em Arch. Rational Mech. Anal.}, 77(1):11--21, 1981.

\bibitem[Bob88]{Bo88}
A.~V. Bobyl{\"e}v.
\newblock The theory of the nonlinear spatially uniform {B}oltzmann equation
  for {M}axwell molecules.
\newblock In {\em Mathematical physics reviews, Vol.\ 7}, volume~7 of {\em
  Soviet Sci. Rev. Sect. C Math. Phys. Rev.}, pages 111--233. Harwood Academic
  Publ., Chur, 1988.

\bibitem[CGT99]{CGT}
E.~A. Carlen, E.~Gabetta, and G.~Toscani.
\newblock Propagation of smoothness and the rate of exponential convergence to
  equilibrium for a spatially homogeneous {M}axwellian gas.
\newblock {\em Comm. Math. Phys.}, 199(3):521--546, 1999.

\bibitem[Des95]{des1}
L.~Desvillettes.
\newblock About the regularizing properties of the non-cut-off {K}ac equation.
\newblock {\em Comm. Math. Phys.}, 168(2):417--440, 1995.

\bibitem[Des03]{des2}
L.~Desvillettes.
\newblock About the use of the {F}ourier transform for the {B}oltzmann
  equation.
\newblock {\em Riv. Mat. Univ. Parma (7)}, 2*:1--99, 2003.
\newblock Summer School on ``Methods and Models of Kinetic Theory'' (M\&MKT
  2002).

\bibitem[DM05]{d-m}
L.~Desvillettes and C.~Mouhot.
\newblock About {$L\sp p$} estimates for the spatially homogeneous {B}oltzmann
  equation.
\newblock {\em Ann. Inst. H. Poincar\'e Anal. Non Lin\'eaire}, 22(2):127--142,
  2005.

\bibitem[Gus86]{gus1}
T.~Gustafsson.
\newblock {$L\sp p$}-estimates for the nonlinear spatially homogeneous
  {B}oltzmann equation.
\newblock {\em Arch. Rational Mech. Anal.}, 92(1):23--57, 1986.

\bibitem[Gus88]{gus2}
T.~Gustafsson.
\newblock Global {$L\sp p$}-properties for the spatially homogeneous
  {B}oltzmann equation.
\newblock {\em Arch. Rational Mech. Anal.}, 103(1):1--38, 1988.

\bibitem[MV04]{mvv}
C.~Mouhot and C.~Villani.
\newblock Regularity theory for the spatially homogeneous boltzmann equation
  with cut-off.
\newblock {\em Arch. Ration. Mech. Anal.}, 173(2):169--212, 2004.

\bibitem[PT96]{PT}
A.~Pulvirenti and G.~Toscani.
\newblock The theory of the nonlinear {B}oltzmann equation for {M}axwell
  molecules in {F}ourier representation.
\newblock {\em Ann. Mat. Pura Appl. (4)}, 171:181--204, 1996.

\bibitem[Rod93]{rod}
L.~Rodino.
\newblock {\em Linear partial differential operators in {G}evrey spaces}.
\newblock World Scientific Publishing Co. Inc., River Edge, NJ, 1993.

\bibitem[TV99]{TV}
G.~Toscani and C.~Villani.
\newblock Probability metrics and uniqueness of the solution to the {B}oltzmann
  equation for a {M}axwell gas.
\newblock {\em J. Statist. Phys.}, 94(3-4):619--637, 1999.

\bibitem[Uka84]{U84}
S.~Ukai.
\newblock Local solutions in {G}evrey classes to the nonlinear {B}oltzmann
  equation without cutoff.
\newblock {\em Japan J. Appl. Math.}, 1(1):141--156, 1984.

\bibitem[Vil02]{vil}
C.~Villani.
\newblock A review of mathematical topics in collisional kinetic theory.
\newblock In {\em Handbook of mathematical fluid dynamics, Vol. I}, pages
  71--305. North-Holland, Amsterdam, 2002.

\end{thebibliography}
\end{document}